\documentclass[a4paper]{amsart}
\usepackage{amsmath,amssymb}
\usepackage[all]{xy}
\usepackage[latin1]{inputenc}
\usepackage{url}
\usepackage{tikz} 
\usetikzlibrary{patterns} 

\usepackage[colorlinks=true]{hyperref}

\newcommand{\ie}{{i.e.}\ }

\newcommand{\ZZ}{{\mathbb Z}} 

\newcommand{\PP}{{\mathbb P}} 
\newcommand{\Gr}{{\mathbb G}r} 
\DeclareMathOperator{\Pic}{Pic} 
\newcommand{\cV}{{\mathcal V}} 
\newcommand{\cO}{{\mathcal O}} 

\newcommand{\W}{W} 

\newcommand{\Dyn}{\Delta}
\newcommand{\LB}{\mathcal{L}} 
\DeclareMathOperator{\SL}{SL} 
\DeclareMathOperator{\GL}{GL} 
\DeclareMathOperator{\Rad}{Rad} 
\newcommand{\Gm}{\mathbb{G}_m} 

\swapnumbers 
\theoremstyle{plain}
\newtheorem*{thm*}{Theorem}

\newtheorem*{lem*}{Lemma}

\theoremstyle{definition}
\newtheorem*{rmk*}{Remark}
\newtheorem*{exa*}{Example}

\title{Trivial Witt groups of flag varieties}
\author{Baptiste Calm{\`e}s and Jean Fasel}

\address{Baptiste Calm\`es, Universit\'e d'Artois, Laboratoire de Math\'ematiques de Lens, France}
\urladdr{http://www.math.uni-bielefeld.de/~bcalmes}

\address{Jean Fasel, Ludwig Maximilian Universit\"at, Munich, Germany} 

\date\today

\keywords{Witt group, Flag variety, Projective homogeneous variety}

\subjclass[2010]{19G99, 11E81, 20G10, 20G15, 14M15, 14M17}

\begin{document}

\begin{abstract}
Let $G$ be a split semi-simple linear algebraic group over a field $k$ of characteristic not $2$. Let $P$ be a parabolic subgroup and let $\LB$ be a line bundle on the projective homogeneous variety $G/P$. We give a simple condition on the class of $\LB$ in $\Pic(G/P)/2$ in terms of Dynkin diagrams implying that the Witt groups $\W^i(G/P,\LB)$ are zero for all integers $i$. In particular, if $B$ is a Borel subgroup, then $\W^i(G/B,\LB)$ is zero unless $\LB$ is trivial in $\Pic(G/B)/2$.
\end{abstract}

\maketitle

The main result of this short note is part of a more general work on Witt groups of split projective homogeneous varieties. However, contrary to the rest of the work, it has a rather quick proof that does not require any heavy machinery. Since other people in the subject have expressed interest and used it in computations, we have decided to write it down as a standalone result. 

The first named author would like to thank Ian Grojnowski for a fruitful discussion.
\medskip

Let $G$ be a split simply connected semi-simple linear algebraic group over a field $k$, $T$ a maximal split torus of $G$ and $B$ a Borel subgroup containing $T$. Let $\Dyn=\Dyn_G$ denote the corresponding Dynkin diagram of $G$, in which the vertices are the simple roots of $G$ with respect to $B$.  Let $\Theta$ be a subset of the vertices of the Dynkin diagram $\Dyn$. The standard parabolic subgroup $P_{\Theta}$ is defined as the subgroup of $G$ generated by $B$ and the root subgroups $U_{-\alpha}$ for all $\alpha \in \Theta$. Any parabolic subgroup of $G$ is conjugate to such a standard parabolic subgroup, so all projective homogeneous varieties under $G$ are of the form $X_\Theta = G/P_{\Theta}$ for some subset $\Theta$ of $\Dyn$. Note that a projective homogeneous variety under a semi-simple group can always be considered as a projective homogeneous variety under the simply connected cover of the group, so the assumption that $G$ is simply connected is harmless for our purposes.
 
To an element $\lambda$ of the weight lattice of the root system of $G$, we associate the line bundle $\LB_\lambda$ over $G/B$ defined as follows. Since $G$ is simply connected, the weight $\lambda$ is a character of $T$ and extends (uniquely up to isomorphism) to a one-dimensional representation $V_\lambda$ of $B$. We define $\LB_\lambda$ as the line bundle whose total space is $G\times_{B} V_\lambda$. The assignment $\lambda \mapsto [\LB_\lambda]$ defines a group isomorphism between the weight lattice and $\Pic(G/B)$. See for example \cite[\S 2.1]{Merkurjev95}. 

To each simple root $\alpha$ corresponds a weight $\omega_\alpha$ characterized by the relation $\langle\omega_\alpha,\beta^\vee \rangle=\delta_{\alpha,\beta}$ for all simple roots $\beta$. These weights are called fundamental weights, and they form a basis of the weight lattice. 
The Picard group of $G/B$ is thus the free abelian group generated by the classes $[\LB_\omega]$ with $\omega$ a fundamental weight. The Picard group of $X_\Theta$ injects in the one of $G/B=X_{\emptyset}$ by pull-back, and its image is the free abelian group generated by the classes of line bundles $\LB_{\omega_\alpha}$ where $\alpha$ is a simple root not in $\Theta$ (see \cite[Prop. 2.3]{Merkurjev95}). We subsequently identify $\Pic(X_{\Theta})$ with its image in $\Pic(G/B)$. This induces a bijection between the subsets of $\Dyn - \Theta$ and the Picard group modulo $2$ of $X_\Theta$, sending a subset $\Lambda$ to the class $[\LB_{\sum_{\alpha \in \Lambda} \omega_\alpha}]$ in $\Pic(X_\Theta)/2$. For a line bundle $\LB$ on $X_\Theta$, we use the notation $\Lambda(\LB)$ for the inverse bijection applied to its class in $\Pic(X_{\Theta})/2$. In other words, a simple root $\alpha$ is in $\Lambda(\LB)$ if and only if $[\LB_{\omega_{\alpha}}]$ appears with an odd multiplicity in the decomposition of $[\LB] \in \Pic(G/B)$.  

We say that a simple root $\alpha \in \Dyn -\Theta$ {\em is not adjacent to} $\Theta$ if no edge of $\Dyn$ connects $\alpha$ and a vertex of $\Theta$, \ie if $\alpha$ is orthogonal to all simple roots $\beta \in \Theta$.

\begin{thm*}
Assume the characteristic of $k$ is not $2$. Let $\LB$ be a line bundle on $X_\Theta$ and let $\Lambda=\Lambda(\LB)$ be the associated subset of $\Dyn -\Theta$. If there is a vertex $\alpha \in\Lambda$ that is not adjacent to $\Theta$, then $\W^i(X_{\Theta},\LB)=0$ for all integers $i$. 
\end{thm*}
\begin{proof}
Let $\Theta'$ be the subset $\Theta \cup \{\alpha\}$. Then, there is a vector bundle $\cV$ of rank $2$ over $X_{\Theta'}$ and an isomorphism $X_\Theta \simeq \PP_{X_{\Theta'}}(\cV)$ of schemes over $X_{\Theta'}$. See the lemma below for a proof.
The theorem is then implied by the following fact: for any vector bundle $\cV$ of rank $2$ over a regular base $X$ and any line bundle $\LB$ on $\PP_X(\cV)$ such that $[\LB]\in \Pic(\PP_X(\cV))/2$ is not in the image of the pull-back from $\Pic(X)/2$, the Witt groups $\W^i(\PP_X(\cV),\LB)$ are zero for all values of $i$ \cite[Thm. 1.3]{Walter03}. 
Obviously, in our case, $[\LB]$ is not in the image of the pull-back from $\Pic(X_{\Theta'})/2$, since $\Pic(X_\Theta)/2=\Pic(X_{\Theta'})/2 \oplus \ZZ/2\ [\LB_{\omega_\alpha}]$ and $\alpha \in \Lambda$ precisely means that the component of $\LB$ on $[\LB_{\omega_\alpha}]$ is nonzero.
\end{proof}
\begin{exa*}
Let $B$ be a Borel subgroup of $G$. Then $\W^i(G/B,\LB)=0$ for any $i$ unless $[\LB] = [\cO] \in \Pic(G/B)/2$. Indeed, $G/B=X_{\emptyset}$ so any $\LB$ such that $\Lambda(\LB) \neq \emptyset$ satisfies the assumption of the theorem. Here are more examples with a group $G$ of type $D_4$. \\

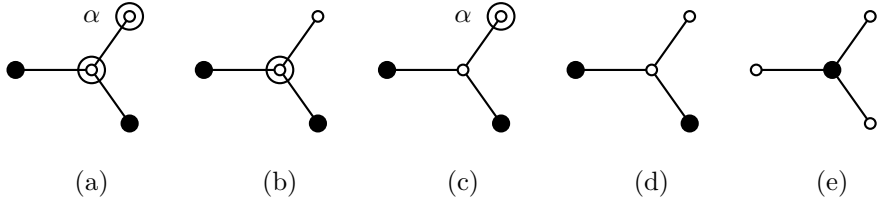
\begin{tikzpicture}[scale=1]
\draw[thick] (0,0) -- (1,0);
\draw[thick] (1,0) -- (1.5,0.71);
\draw[thick] (1,0) -- (1.5,-0.71);
\foreach \x / \y in {0/0,1/0,1.5/0.71,1.5/-0.71}
\filldraw[thick,fill=white] (\x,\y) circle (2pt);
\foreach \x / \y in {0/0,1.5/-0.71}
\filldraw[thick,fill=black] (\x,\y) circle (3pt);
\foreach \x / \y in {1/0,1.5/0.71}
\draw[thick,] (\x,\y) circle (5pt);
\draw (1,0.71) node{$\alpha$}; 
\draw (1,-1.5) node{(a)}; 
\end{tikzpicture}
\hfil
\begin{tikzpicture}[scale=1]
\draw[thick] (0,0) -- (1,0);
\draw[thick] (1,0) -- (1.5,0.71);
\draw[thick] (1,0) -- (1.5,-0.71);
\foreach \x / \y in {0/0,1/0,1.5/0.71,1.5/-0.71}
\filldraw[thick,fill=white] (\x,\y) circle (2pt);
\foreach \x / \y in {0/0,1.5/-0.71}
\filldraw[thick,fill=black] (\x,\y) circle (3pt);
\foreach \x / \y in {1/0}
\draw[thick,] (\x,\y) circle (5pt);
\draw (1,-1.5) node{(b)}; 
\end{tikzpicture}
\hfil
\begin{tikzpicture}[scale=1]
\draw[thick] (0,0) -- (1,0);
\draw[thick] (1,0) -- (1.5,0.71);
\draw[thick] (1,0) -- (1.5,-0.71);
\foreach \x / \y in {0/0,1/0,1.5/0.71,1.5/-0.71}
\filldraw[thick,fill=white] (\x,\y) circle (2pt);
\foreach \x / \y in {0/0,1.5/-0.71}
\filldraw[thick,fill=black] (\x,\y) circle (3pt);
\foreach \x / \y in {1.5/0.71}
\draw[thick,] (\x,\y) circle (5pt);
\draw (1,0.71) node{$\alpha$}; 
\draw (1,-1.5) node{(c)}; 
\end{tikzpicture}
\hfil
\begin{tikzpicture}[scale=1]
\draw[thick] (0,0) -- (1,0);
\draw[thick] (1,0) -- (1.5,0.71);
\draw[thick] (1,0) -- (1.5,-0.71);
\foreach \x / \y in {0/0,1/0,1.5/0.71,1.5/-0.71}
\filldraw[thick,fill=white] (\x,\y) circle (2pt);
\foreach \x / \y in {0/0,1.5/-0.71}
\filldraw[thick,fill=black] (\x,\y) circle (3pt);
\draw (1,-1.5) node{(d)}; 
\end{tikzpicture}
\hfil
\begin{tikzpicture}[scale=1]
\draw[thick] (0,0) -- (1,0);
\draw[thick] (1,0) -- (1.5,0.71);
\draw[thick] (1,0) -- (1.5,-0.71);
\foreach \x / \y in {0/0,1/0,1.5/0.71,1.5/-0.71}
\filldraw[thick,fill=white] (\x,\y) circle (2pt);
\foreach \x / \y in {1/0}
\filldraw[thick,fill=black] (\x,\y) circle (3pt);
\draw (1,-1.5) node{(e)}; 
\end{tikzpicture}
\def\thetapoint{\tikz[baseline=-2pt]\filldraw[thick,fill=black] (0,0) circle (3pt);}
\def\lambdapoint{\tikz[baseline=-2pt]\draw[thick] (0,0) circle (5pt) circle (2pt);}

\noindent The vertices in $\Theta$ are marked \protect\thetapoint, while the vertices in $\Lambda(\LB)$ are marked \protect\lambdapoint. Cases (a) and (c) satisfy the assumption of the theorem (with the given $\alpha$), while cases (b) and (d) do not. With $\Theta$ as in case (e), no $\Lambda(\LB)$ can satisfy the assumption of the theorem.
\end{exa*}

\begin{rmk*}
The converse of the theorem is not true: it can happen that a line bundle $\LB$ is such that $\W^i(X_\Theta,\LB)=0$ for all values of $i$, while all $\alpha \in \Lambda(\LB)$ are adjacent to $\Theta$. Examples are provided by quotients of a group of type $A_n$ by a maximal proper parabolic subgroup: all vertices are in $\Theta$ but the $d$-th one. This yields a Grassmann variety $\Gr(d,d+e)$ with $n=d+e-1$. When $d$ and $e$ are odd, and when $\LB$ is a generator of the Picard group of $\Gr(d,d+e)$, then $\W^i(\Gr(d,d+e),\LB)=0$ for all values of $i$ (see \cite{Balmer07}).
\end{rmk*}

Since we could not find a reference in the literature for the following well-known fact, we include a proof (actually valid over any base, not necessarily a field).
\begin{lem*}
If $\alpha$ is a simple root orthogonal to all simple roots in $\Theta$, then the natural projection $X_\Theta \to X_{\Theta \cup \{\alpha\}}$ is isomorphic to the projective bundle associated to a vector bundle of rank $2$ over $X_{\Theta \cup \{\alpha\}}$. 
\end{lem*}
\begin{proof}
Recall that $G$ is simply connected. Let $P=P_\Theta$ and $Q=P_{\Theta \cup \{\alpha\}}$. Let $L$ be the Levi subgroup of $Q$ containing the maximal torus $T$. Its derived subgroup $\left[L,L\right]$ is semi-simple, simply connected and has Dynkin type given by the full subgraph of $\Dyn$ whose vertices are in $\Theta\cup \{\alpha\}$ \cite[Prop. 6.2.7]{SGA3-3}. 
Thus, $[L,L]\simeq \SL_2 \times H$ with $H$ a semi-simple simply connected subgroup of $P$. The radical, $\Rad(L)$, is a central torus in $L$ and the multiplication map $\left[L,L\right] \times \Rad(L) \to L$ is surjective \cite[6.2.3]{SGA3-3}. Its kernel $\nu$ is isomorphic to $\left[L,L\right] \cap \Rad(L)$, and is therefore a finite subgroup of the center $\mu_2 \times \mu$ of $\left[L,L\right]$, where $\mu$ is a finite group of multiplicative type. There are two cases: (1) the kernel $\nu$ is included in $\mu$ or (2) it surjects to $\mu_2$ through the projection map $\mu_2 \times \mu \to \mu_2$. In the first case, the projection $\SL_2 \times H \times \Rad(L) \to \SL_2$ factors as a (surjective) map $L \to \SL_2$. In the second case, let us show that $\Rad(L)$ has a surjection to $\Gm$ such that the following diagram on the left commutes.
$$\xymatrix{
\nu \ar@{^{(}->}[r] \ar@{->>}[d] & \Rad(L) \ar@{->>}[d] \\ 
\mu_2 \ar@{^{(}->}[r] & \Gm
} \hspace{10ex}
\xymatrix{
N & \ZZ^{r} \ar@{->>}[l] \\
\ZZ/2 \ar@{^{(}->}[u] & \ZZ \ar@{->>}[l] \ar@{-->}[u]
}$$
Taking Cartier duals, it amounts to checking the existence of a dashed map such that the right diagram commutes: take the image of $1 \in \ZZ/2$ in $N$, lift it to $\ZZ^r$ and send $1 \in \ZZ$ to this lift. This implies that the induced surjection $\SL_2 \times H \times \Rad(L) \to \SL_2 \times \Gm$ factors as a surjection $L \to (\SL_2 \times \Gm)/\mu_2 = \GL_2$. So, in cases (1) and (2), taking the standard representation $V$ of $\SL_2$ (resp. $\GL_2$), and using that $L$ is the quotient of $Q$ by its unipotent radical, we obtain a surjective map $Q \to \SL(V)$ (resp. $\GL(V)$); let $K$ be its kernel. Both $\Rad(L) \subset T \subset P$ and $H \subset P$, so $K\subset P$. Furthermore, $P/K$ is a Borel subgroup of $Q/K$ since $Q/P=L/(P\cap L) \cong \PP^1$. Let us consider the vector bundle $G \times_Q V \to G/Q$. Its associated projective space is $G \times_Q (Q/K)/(P/K) = G\times_Q Q/P =G/P$, as expected.
\end{proof}

\bibliographystyle{alpha}
\bibliography{projhomnullwitt}

\begin{thebibliography}{Dem64}

\bibitem[BC07]{Balmer07}
P.~Balmer and B.~Calm\`es.
\newblock Witt groups of {G}rassmann varieties.
\newblock \href{http://arxiv.org/abs/0807.3296}{arXiv:0807.3296}, 2007.

\bibitem[Dem64]{SGA3-3}
M.~Demazure.
\newblock {\em Sch\'emas en groupes. {III}: {S}tructure des sch\'emas en
  groupes r\'eductifs, Exp. XXII}.
\newblock S\'eminaire de G\'eom\'etrie Alg\'ebrique du Bois Marie 1962/64 (SGA
  3). Dirig\'e par M. Demazure et A. Grothendieck. Lecture Notes in
  Mathematics, Vol. 153. Springer-Verlag, Berlin, 1962/1964.

\bibitem[MT95]{Merkurjev95}
A.~S. Merkurjev and J.-P. Tignol.
\newblock The multipliers of similitudes and the {B}rauer group of homogeneous
  varieties.
\newblock {\em J. Reine Angew. Math.}, 461:13--47, 1995.

\bibitem[Wal03]{Walter03}
C.~Walter.
\newblock Grothendieck-{W}itt groups of projective bundles.
\newblock preprint, \url{http://www.math.uiuc.edu/K-theory/0644/}, 2003.

\end{thebibliography}

\end{document}